 \documentclass[12pt]{article}
\usepackage[dvips]{epsfig}
\usepackage{amsthm,amsmath,amssymb}
\usepackage{latexsym}
\usepackage{placeins}
\usepackage{dblfloatfix}
\usepackage{afterpage}
\usepackage{color}
\usepackage{hyperref}
\hypersetup{
     colorlinks=true,
     linkcolor=cyan,
     filecolor=cyan,
     citecolor = black,      
     urlcolor=cyan,
     }
\definecolor{TWhite}{rgb}{1,1,1}
\definecolor{TBlack}{rgb}{0,0,0}
\definecolor{TBB}{rgb}{0,0,0.4}
\definecolor{LBlack}{rgb}{0.70,0.67,0.70}
\definecolor{Gray}{rgb}{0.55,0.52,0.55}
\definecolor{Comment}{rgb}{1,1,1}
\definecolor{TBlue}{rgb}{0,0,1}
\definecolor{LBlue}{rgb}{0.35,0.60,0.85}
\definecolor{TGreen}{rgb}{0,0.40,0.10}
\definecolor{LGreen}{rgb}{0.0,0.70,0.40}
\definecolor{Yellow}{rgb}{0.55,0.55,0}
\definecolor{TRed}{rgb}{1,0,0}
\definecolor{LRed}{rgb}{1.00,0.4,0.4}
\definecolor{TViolet}{rgb}{0.65,0,0.65}
\definecolor{LViolet}{rgb}{0.95,0.55,0.95}
\definecolor{TBrown}{rgb}{0.45,0.15,0}
\definecolor{LBrown}{rgb}{0.80,0.50,0}

\definecolor{blue}{rgb}{0,0,0.9}
\definecolor{red}{rgb}{0.9,0,0}
\definecolor{green}{rgb}{0,0.50,0.10}
\definecolor{violet}{rgb}{0.5804,0.0000,0.8275}

\def\@themcountersep{}

\newtheorem{THEO}{Theorem}[section]
\newtheorem{ALGo}[THEO]{Algorithm}
\newtheorem{CONJ}[THEO]{Conjecture}
\newtheorem{COND}[THEO]{Condition}
\newtheorem{ASSUMP}[THEO]{Assumption}
\newtheorem{CORO}[THEO]{Corollary}
\newtheorem{DEFI}[THEO]{Definition}
\newtheorem{EXAMP}[THEO]{Example}
\newtheorem{FACT}[THEO]{Fact}
\newtheorem{HYPO}[THEO]{Hypothesis}
\newtheorem{LEMM}[THEO]{Lemma}
\newtheorem{PROB}[THEO]{Problem}
\newtheorem{PROP}[THEO]{Proposition}
\newtheorem{REMA}[THEO]{Remark}
\newcommand{\theo}{\begin{THEO}}
\newcommand{\algo}{\begin{ALGo} \rm}
\newcommand{\cond}{\begin{COND} \rm}
\newcommand{\assump}{\begin{ASSUMP} \rm}
\newcommand{\conj}{\begin{CONJ}}
\newcommand{\coro}{\begin{CORO}}
\newcommand{\defi}{\begin{DEFI} \rm}
\newcommand{\examp}{\begin{EXAMP} \rm}
\newcommand{\fact}{\begin{FACT}}
\newcommand{\hypo}{\begin{HYPO} \rm}
\newcommand{\lemm}{\begin{LEMM}}
\newcommand{\prob}{\begin{PROB} \rm}
\newcommand{\prop}{\begin{PROP}}
\newcommand{\rema}{\begin{REMA} \rm}
\newcommand{\etheo}{\end{THEO}}
\newcommand{\ealgo}{\end{ALGo}}
\newcommand{\econd}{\end{COND}}
\newcommand{\eassump}{\end{ASSUMP}}
\newcommand{\econj}{\end{CONJ}}
\newcommand{\ecoro}{\end{CORO}}
\newcommand{\edefi}{\end{DEFI}}
\newcommand{\eexamp}{\end{EXAMP}}
\newcommand{\efact}{\end{FACT}}
\newcommand{\ehypo}{\end{HYPO}}
\newcommand{\elemm}{\end{LEMM}}
\newcommand{\eprob}{\end{PROB}}
\newcommand{\eprop}{\end{PROP}}
\newcommand{\erema}{\end{REMA}}

\def\0{\mbox{\bf 0}}
\def\1{\mbox{\bf 1}}
\def\2{\mbox{\bf 2}}
\def\3{\mbox{\bf 3}}
\def\4{\mbox{\bf 4}}
\def\5{\mbox{\bf 5}}
\def\6{\mbox{\bf 6}}
\def\7{\mbox{\bf 7}}
\def\8{\mbox{\bf 8}}
\def\9{\mbox{\bf 9}}

\def\b{\mbox{\boldmath $b$}}

\def\u{\mbox{\boldmath $u$}}
\def\v{\mbox{\boldmath $v$}}

\def\x{\mbox{\boldmath $x$}}
\def\y{\mbox{\boldmath $y$}}

\def\A{\mbox{\boldmath $A$}}

\def\F{\mbox{\boldmath $F$}}

\def\H{\mbox{\boldmath $H$}}
\def\I{\mbox{\boldmath $I$}}

\def\O{\mbox{\boldmath $O$}}

\def\Q{\mbox{\boldmath $Q$}}

\def\S{\mbox{\boldmath $S$}}

\def\X{\mbox{\boldmath $X$}}

\def\Z{\mbox{\boldmath $Z$}}

\def\Real{\mbox{$\mathbb{R}$}}

\def\coneN{\mbox{$\mathbb{N}$}}

\def\SymMat{\mbox{$\mathbb{S}$}}

\let\chige\c
\let\dotlessi\i

\begin{document}

%\title{ \Large 
%Conic optimization problems with no duality gap
%} 

%\title{A note on 
%second-order sufficient conditions for the global optimality in QCQOPs and the exactness of 
%their SDP relaxations
%}
\title{
Equivalent Sufficient Conditions for  Global Optimality of Quadratically Constrained 
Quadratic Program}

%\author{
%\normalsize 
%Sunyoung Kim\thanks{Department of Mathematics, Ewha W. University, 52 Ewhayeodae-gil, Sudaemoon-gu, Seoul 03760, Korea 
%			({\tt skim@ewha.ac.kr}). 
%			 The research was supported  by   NRF 2021-R1A2C1003810.}, \and \normalsize
%Masakazu Kojima\thanks{Department of Industrial and Systems Engineering,
%	Chuo University, Tokyo 192-0393, Japan ({\tt kojima@is.titech.ac.jp}).
%	This research was supported by Grant-in-Aid for Scientific Research (A) 19H00808 
% 	}
%}

\author{
\and \normalsize    Sunyoung Kim\footnotemark[1]
\and \normalsize  Masakazu Kojima\footnotemark[2]}
        
\date{\normalsize\today}
\maketitle 

\renewcommand{\thefootnote}{\fnsymbol{footnote}}
\footnotetext[1]{ Corresponding author.}
\footnotetext[1]{Department of Mathematics, Ewha W. University, 52 Ewhayeodae-gil, Sudaemoon-gu,
	Seoul 	03760, Korea  ({\tt skim@ewha.ac.kr}). This work was supported
        by  NRF 2021-R1A2C1003810.}
\footnotetext[2]{
        Department of Industrial and Systems Engineering,
	Chuo University, Tokyo 192-0393, Japan ({\tt kojima@is.titech.ac.jp}).
	This research was supported by Grant-in-Aid for Scientific Research (A) 19H00808.}

% \date{\normalsize Revised, \today}
%\date{\today}

\maketitle 
%\vspace*{-0.4cm}

%\input abstract.tex
%!TEX root = ./main.tex
\begin{abstract}
\noindent
We study the equivalence of several well-known sufficient optimality conditions 
for a general quadratically constrained quadratic program (QCQP). 
The conditions are classified in two categories. The first one is for determining 
an optimal solution and the second one is for finding an optimal value.
The first category of conditions includes 
the existence of a saddle point of the Lagrangian function and the existence of a rank-1 
optimal solution of the primal SDP relaxation of QCQP. The second category   
includes $\eta_p = \zeta$, $\eta_d = \zeta$, and $\varphi = \zeta$, where $\zeta$, $\eta_p$, $\eta_d$, and 
$\varphi$ denote the optimal values of QCQP, the dual SDP relaxation, the primal SDP relaxation 
and the Lagrangian dual, respectively.
We show the equivalence of these conditions with or without the existence of an optimal solution of QCQP 
and/or the Slater constraint qualification for the primal SDP relaxation.
The results on the conditions are also extended to the doubly nonnegative relaxation
of equality constrained QCQP in nonnegative variables.
\end{abstract}
%
%\vspace{0.5cm}
%
%\noindent
%
\bigskip

%\vspace{0.5cm}
%

%\input sect1.tex
%!TEX root = ./main.tex
\section{Introduction}

We consider a {\em quadratically constrained quadratic program} (QCQP). 
\begin{eqnarray}
\zeta & = &  \inf \left\{q_0(\u): q_k(\u) \leq 0 \ (k \in M)
\right\}, \label{eq:QCQP00} 
\end{eqnarray}
where 
\begin{eqnarray*}
%\left. 
%\begin{array}{l}
& & M = \{1,\ldots,m\}, \\ [2pt]
& & q_k(u) = \u^T\A_k\u + 2\b_k^T\u + c_k \ \mbox{for every } \u \in \Real^n, \\ [2pt]
& & \A_k \in  \SymMat^n \ \mbox{(the linear space of $n \times n$ matrix)}, \\ [2pt]
& & \b_k \in \Real^n \ \mbox{(the $n$-dim. Euclidean space of column vector)}, \\ [2pt]
& & c_k  \in  \Real, \ c_0 \ = \ 0,
%\end{array} 
%\right\} 
\end{eqnarray*}
$(0 \leq k \leq m)$. 
Let $F$ denote the feasible 
region of QCQP~\eqref{eq:QCQP00}; 
$F = \{\u \in \Real^n:q_k(\u)\leq 0 \ (k \in M)\}$. 
We call each $\u \in F$ {\em a feasible solution} of QCQP~\eqref{eq:QCQP00}, 
$\u \in F$ satisfying $q_0(\u) \leq 
q_0(\u')$ for all $\u' \in \{ \u' \in F : \|\u' -\u\| <\epsilon\}$ and 
some $\epsilon > 0$ {\em a local 
minimizer} of QCQP~\eqref{eq:QCQP}, and $\u \in F$ satisfying $q_0(\u) \leq q_0(\u')$ for 
all $\u' \in F$ 
{\em a global minimizer} of QCQP~\eqref{eq:QCQP00}. Obviously, a global minimizer 
is a local minimizer. 

QCQP is one of the most fundamental nonconvex  
optimization problems that include various important NP-hard 
 problems, notably,  
max cut problems \cite{GEOMANS95}, maximum stable set problems 
\cite{DEKLERK2002}, graph partitioning problems \cite{POVH2007}, 
and quadratic assignment problems \cite{LOIOLA2007}. 
It is also known \cite{MEVISSEN2010} that any polynomial optimization problem can be converted into a QCQP. 
For NP-hard QCQP, finding the exact optimal solution  or the exact optimal value  is an 
important issue. 

Our focus here is on conditions which characterize  global optimality 
for  QCQP~\eqref{eq:QCQP00};  more precisely, three conditions 
(Conditions (A), (B) and (C)) for a feasible solution 
$\u$ of~\eqref{eq:QCQP00} to be a global minimizer of~\eqref{eq:QCQP00}, and 
three conditions (Conditions (D), (E) and (F)) 
on some lower bound of the optimal value $\zeta$ of QCQP to be tight. 
Specifically, the main purpose of this paper is to 
clarify and understand their relations by showing that each of them is equivalent to all the others with 
or without additional moderate assumptions. To describe these conditions, we use \vspace{-2mm}
\begin{itemize}
\item the Lagrangian function $L$ in the variable vector $\u \in \Real^n$ and the 
multiplier vector $\y\in \Real^m$,\vspace{-2mm} 
\item the  primal semidefinite (SDP) relaxation~\eqref{eq:PSDP} of QCQP~\eqref{eq:QCQP00} with the optimal value $\eta_p$
in an $(1+n) \times (1+n)$ symmetric variable matrix $\X$,\vspace{-2mm} 
\item  the dual SDP relaxation~\eqref{eq:DSDP} of~\eqref{eq:QCQP00} 
 with the optimal value $\eta_d$ in a variable vector $(\y,s) \in \Real^{m+1}$ . \vspace{-2mm} 
\end{itemize}

Let $G = \{ \y \in \Real^m: y_k \geq 0 \ (k\in M) \}$ and 
$L:\Real^{n+m} \rightarrow \Real$ 
be the (standard) {\em Lagrangian function} for QCQP~\eqref{eq:QCQP00} defined by 
\begin{eqnarray*}
L(\u,\y) & = & q_0(\u) + \sum_{k\in M} y_k q_k(\u) \ \mbox{for every }
 (\u,\y) \in \Real^{n+m}. 
\end{eqnarray*} 
Condition (A) is described through the {\em saddle point problem}: 
Find a $(\u,\y) \in \Real^{n}\times G$ such that 
\begin{eqnarray}
\sup_{\y' \in G}L(\u,\y') = L(\u,\y) = \inf_{\u' \in \Real^n} L(\u',\y). 
\label{eq:saddlePoint}
\end{eqnarray}
This problem was introduced in the book~\cite{MANGASARIAN1969} as a sufficient 
condition for $\u$ to be a global minimizer 
of a more general optimization problem where 
$q_k : \Real^n \rightarrow \Real$ $(0 \leq k \leq m)$ are allowed to be continuous functions:  
if $(\u,\y)$ is a solution of~\eqref{eq:saddlePoint}, then $\u$ is a minimizer 
of~\eqref{eq:QCQP00} 
\cite[5.3.1]{MANGASARIAN1969}. 
We also consider the {\em Lagrangian dual}: 
\begin{eqnarray}
\varphi & = & \sup_{\y\in G} \inf\left\{L(\u,\y) : \u \in \Real^n \right\}. \label{eq:LagrangianDual}
\end{eqnarray}
It is well-known that $\varphi$ provides a lower bound of $\zeta$. 
We let Condition (F) be `$\varphi = \zeta$'.

All the other conditions are described through the primal SDP relaxation~\eqref{eq:PSDP}  
and the dual SDP relaxation~\eqref{eq:DSDP} of QCQP~\eqref{eq:QCQP00} 
\cite{BAO2011,FUJIE1997,SHOR1987,SHOR1990}.
In general, their optimal values $\eta_p$ and $\eta_d$ satisfy $\eta_d \leq \eta_p \leq \zeta$. 
Conditions (D) and (E) are `$\eta_p = \zeta$' and `$\eta_d = \zeta$', respectively. 
If the primal SDP relaxation of QCQP~\eqref{eq:QCQP00} can provide  
a minimizer $\u$ of QCQP~\eqref{eq:QCQP00}, then we call the SDP relaxation {\em exact}. 
Classes of QCQPs whose primal SDP (and/or second order cone programming (SOCP)) relaxations 
are exact have been studied 
extensively in \cite{AZUMA2022a,AZUMA2022,KIM2003,SOJOUDI2014,Wang2021tightness,ZHANG2000}, where 
the minimizer $\u$ can be derived from a rank-$1$ optimal solution $\X$ of the primal SDP relaxation 
with the form $\X = \begin{pmatrix}1 \\ \u \end{pmatrix}\begin{pmatrix}1 \\ \u \end{pmatrix}^T$. 
% It was shown in \cite{DUR2017} that the Slater condition and strong duality is a generic 
% property of conic problems.
Each QCQP in those classes has been identified by its data matrices $\A_k$ and vectors $\b_k$ $(0 \leq k \leq m)$ 
that satisfy a certain structured sparsity such as 
tridiagonal, forest and bipartite 
and/or a certain sign-definiteness property. 
In addition, strong duality was assumed in  \cite{AZUMA2022a,AZUMA2022}.
In \cite{Jeyakumar2014}, the exact SDP relaxation of an extended trust-region type 
QCQP was studied under a certain dimension condition.
In \cite{BURER2020}, a general QCQP with no particular structure was transformed 
to a diagonal QCQP whose primal SDP relaxation is exact. 
Condition (C) is  `the primal SDP relaxation~\eqref{eq:PSDP} is exact'.

The relations among Conditions (A) through (F) shown in this paper are summarized as follows:
\begin{eqnarray}
\left.
\begin{array}{lcl}
& & \hspace{-2mm} \mbox{(A)  $\exists$ a saddle point $(\u,\y)$ of the Lagrangian function $L$}. \\[2pt]
& & \hspace{-1mm} 
\Updownarrow \\[2pt]
& & \hspace{-2mm} \mbox{(B) the primal SDP~\eqref{eq:PSDP} is exact}, \\[2pt] 
& & \hspace{6mm} \mbox{$\exists$ an optimal solution $(\y,s)$ of the dual SDP~\eqref{eq:DSDP} and}\\[2pt]
& & \hspace{6mm} \mbox{$\eta_d = \eta_p$ (the strong duality)}. \\[2pt]
& & % \hspace{2mm} 
\Downarrow \hspace{1mm} \Uparrow \mbox{if (a) the 
Slater constraint qualification~\eqref{eq:SlatrCondition} holds.}\\[2pt]
& & \hspace{-2mm} \mbox{(C) the primal SDP~\eqref{eq:PSDP} is exact.}\\[2pt]
& & % \hspace{2mm} 
\Downarrow \hspace{2mm} \Uparrow \mbox{if (b) QCQP~\eqref{eq:QCQP00} has a minimizer $\u$ when $\zeta$ is finite.}\\[2pt]
& & \hspace{-2mm} \mbox{(D) $-\infty < \eta_p = \zeta  < \infty$.}\\[2pt]
& \mbox{if (B) holds} & % \hspace{-1mm} 
\Downarrow \hspace{3mm} \Uparrow \\[2pt]
& & \hspace{-2mm} \mbox{(E) $-\infty < \eta_d = \zeta < \infty$} \Longleftrightarrow 
\mbox{(F) $-\infty < \varphi = \zeta < \infty$.} 
\end{array}
\right\}
\label{eq:relations}
\end{eqnarray}
Apparently, equivalent (A) and (B) are the strongest conditions, and Condition (D) the weakest. 
We see that all the conditions are equivalent if (a) the 
Slater constraint qualification~\eqref{eq:SlatrCondition} holds and if (b) QCQP~\eqref{eq:QCQP00} has a minimizer 
$\u$ when $\zeta$ is finite. 
Since the assumptions (a) and (b) are regarded as 
to avoid special degenerate cases, it can be  approximately said 
 that all Conditions (A) through (F) are equivalent except 
special degenerate cases. 
%%%
In fact, it was shown in \cite{DUR2017} that the Slater condition 
is a generic property of conic optimization problems. Also, 
if the feasible region of QCQP~\eqref{eq:QCQP00} is bounded, then (b) holds.
%%%
We present some examples for  such exceptional degenerate cases. \vspace{-2mm}
\begin{itemize}
\item (C) $\not \Rightarrow$ (B) (Example~\ref{examp:QCQP41}), \vspace{-2mm}
\item (D) $\not \Rightarrow$ (C) (Example~\ref{examp:QCQP43}), \vspace{-2mm}
\item (C), (D) and QCQP~\eqref{eq:QCQP00} has a minimizer $\not 
\Rightarrow$ (E) (Example~\ref{examp:QCQP44}). \vspace{-2mm}
\end{itemize}
%}

\bigskip
\noindent
{\bf Some related works.} 
In general, the class of QCQPs whose SDP relaxation is exact is limited as mentioned above.
Sufficient global optimality conditions on QCQP via the SDP relaxation are not 
strong enough to cover 
the entire class of general QCQPs. 
Some stronger convex conic programming relaxations have been proposed for other classes of QCQPs. They provide 
a lower bound for the optimal value of QCQP, so they serve as a sufficient global optimality condition 
for general QCQPs. 
%One 
A stronger convex relaxation is the completely positive programming cone (CPP) relaxation. 
It is known that CPP relaxation is 
exact for a class of QCQPs with linear and complementarity constraints in nonnegative 
continuous and/or binary variables 
\cite{BURER2009,DUR2021,KIM2020}. CPP relaxation is, however, mainly of theoretical interest 
since the CPP relaxation problem is NP hard. The doubly nonnegative (DNN) relaxation 
\cite{KIM2013,KIM2021,YOSHISE2010} is a numerically implementable 
relaxation of the CPP relaxation, which is at least as strong as the SDP relaxation. 
It was shown in \cite{KIM2020a} that the DNN relaxation is exact for a class of QCQPs with block-clique structure. 
In their paper \cite{LU2011}, Lu et al. proposed an equivalent reformulation of a general QCQP, 
which may be regarded as a (strengthening) modification of the CPP relaxation. They further relaxed their  
modified relaxation to a numerically implementable one which aims
 to compute a global minimizer. 

\bigskip

\noindent
{\bf Contribution.} The main contribution is  %~\eqref{eq:relations} 
to show the equivalence or inclusion  relations  
among Conditions (A) through (F) on  global optimality of QCQP illustrated in~\eqref{eq:relations}. 
While some part of the relations may appear in a scattered manner in the literature, the comprehensive  
relations among the conditions have not been presented.
%All of the relations may be well-known, but they have been presented in diversified literature. 
 %the reader can grasp 
With~\eqref{eq:relations}, the entire
equivalence relations with or without moderate additional assumptions  can be clearly understood. 
Moreover, Examples 4.1 through 4.4 show  some exceptional cases where the equivalence relation does not hold. 
%cases.  
% }

\bigskip

This paper is organized as follows: 
Some notation and symbols used throughout this paper are listed in Section 2.1. 
We present a global optimality condition via the saddle point problem 
for a general nonlinear program in Section 2.2, which corresponds to Condition (A), 
and specialize it to a global optimality condition, Condition (A') for QCQP~\eqref{eq:QCQP00} in Section 2.3. 
In Section 2.4, we introduce the primal SDP relaxation~\eqref{eq:PSDP} and  the  dual SDP relaxation~\eqref{eq:DSDP} 
of QCQP~\eqref{eq:QCQP00},  and present a well-known 
sufficient optimality condition, the Karush-Kuhn-Tucker (KKT) condition. We then combine the 
KKT condition with Condition (C) `the  primal SDP~\eqref{eq:PSDP} is exact' for Condition (B'), which is equivalent to 
Condition (B). 
Section 3 is devoted to proofs of all relations in~\eqref{eq:relations}. 
Four examples, Examples 4.1 through 4.4 are presented in Section 4. 
Section 5 extends Conditions (A) through (F) to an equality constrained QCQP in nonnegative variables 
with DNN relaxation.
We give concluding remarks in Section 6.

%
%\input sect2.tex
%!TEX root = ./main.tex
\section{Preliminaries}

\subsection{Notation and symbols}

Let $\Real$ denote the set of real numbers, $\Real^n$ 
the $n$-dimensional Euclidean space of column vectors 
$\u = (u_1,\ldots,u_n)$ with elements $u_i \in \Real$ $(1 \leq i \leq n)$, and 
$\SymMat^n$ the linear space of $n \times n$ symmetric matrices 
$\A = [A_{ij}]$ with elements $A_{ij} \in \Real$ $(1 \leq i,j \leq n)$. 
The row vector  $\u^T$ stands for the transposed vector of $\u$ for every $\u \in \Real^n$. 
We assume that if $\x \in \Real^{1+n}$ and/or $\X \in \SymMat^{1+n}$,  then their column 
and row indices run from $0$ through $n$, {\it i.e.}, $\x = (x_0,x_1,\ldots,x_n)$ 
and the elements of $\X$ are $X_{ij}$ $(0 \leq i,j \leq n)$. 
For $\Q, \X \in \SymMat^{1+n}$, their inner product is written as
$\Q \bullet \X = \sum_{i=0}^{n}\sum_{j=0}^{n} Q_{ij}X_{ij}$. 
Let 
\begin{eqnarray*}
\SymMat^{\ell}_+&=&\mbox{the cone of positive semidefinite matrices in $\SymMat^{\ell}$},\\
\SymMat^{\ell}_{++}&=&\mbox{the cone of positive matrices in $\SymMat^{\ell}$},  
\end{eqnarray*}
where $\ell = n$ or $1+n$. The zero vector and zero matrix are denoted by $\0$, the $n$-dimensional column vector 
with all elements $0$, and $\O$, the $n \times n$ matrix with all elements $0$, respectively. For each twice 
continuously differentiable function $f: \Real^n \rightarrow \Real$, $\nabla_{u} f(\u)$ 
denotes the gradient row vector of $f$ with elements $\displaystyle \frac{d f(\u)}{d u_i}$ 
$(1 \leq i \leq n)$, and $\nabla_{uu}f(\u)$ 
the $n \times n$ Hessian matrix of $f$ with elements 
$\displaystyle \frac{d^2 f(\u)}{d u_i d u_j}$ $(1 \leq i,j\leq n)$. 

\subsection{Global optimality via the saddle point problem in general nonlinear programs}

Throughout 
Section 2.2, we assume that $q_k : \Real^n \rightarrow \Real$
$(0 \leq k \leq m)$ are twice continuously differentiable functions, but not necessarily quadratic. 
Given $\y \in \Real^m$, 
we denote the $n$-dimensional gradient row vector and the $n \times n$ Hessian matrix 
of $L(\cdot,\y) : \Real^n \rightarrow \Real$ evaluated at $\u \in \Real^n$ by 
$\nabla_{u}L(\u,\y)$ and $\nabla_{uu}L(\u,\y)$, respectively; 
\begin{eqnarray*}
\nabla_{u}L(\u,\y) & = & \nabla_{u}q_0(\u) + \sum_{k\in M} y_k \nabla_{u}q_k(\u), \\
\nabla_{uu}L(\u,\y) & = & \nabla_{uu}q_0(\u) + \sum_{k\in M} y_k \nabla_{uu}q_k(\u).
\end{eqnarray*}

We note that the right equality of the saddle point problem~\eqref{eq:saddlePoint} 
corresponds to the {\em Lagrangian relaxation problem} 
\begin{eqnarray*}
\varphi(\y) = \inf \{ L(\u,\y) : \u \in \Real^n \}. 
\end{eqnarray*}
On the left side of~\eqref{eq:saddlePoint}, we observe that 
$\sup \left\{ L(\u,\y') : \y \in G \right\} = \infty$ if $\u \not\in F$, and that  
$L(\u,\y) < L(\u,\0)$ if $\u \in F$ and $y_kq_k(\u) < 0$ for some $k \in M$. 
Hence, if the left side of~\eqref{eq:saddlePoint} holds, then  
\begin{eqnarray}
\u \in F, \ y_kq_k(\u) = 0 \ (k \in M), 
\label{eq:complemetarity} 
\end{eqnarray}
which implies that $L(\u,\y) = q_0(\u)$. 
It is straightforward to verify  that the converse is true; hence they are equivalent. 
Therefore, we obtain that $\u \in F$ and $\zeta \leq q_0(\u)=L(\u,\y) = \varphi(\y) \leq \varphi \leq \zeta$ or equivalently $\varphi(\y) =\varphi = \zeta = q_0(\u)=L(\u,\y)$.   
By the discussion above, we know that 
Condition (A) is sufficient for $\u$ to be a global minimizer of~\eqref{eq:QCQP00},  
and that (A) holds 
if and only if $L(\u,\y) = \varphi(\y)$ (or $L(\u,\y) = \varphi$) and~\eqref{eq:complemetarity} hold.

\subsection{Global optimality in QCQP~\eqref{eq:QCQP00}}

We apply Condition (A) specifically to QCQP~\eqref{eq:QCQP00} with 
quadratic $q_k$ $(0 \leq k \leq m)$. In this case, we see that 
\begin{eqnarray*}
L(\u,\y) & = & \u^T\A_0\u + 2 \b_0^T\u + \sum_{k\in M} y_k (\u^T\A_k\u + 2 \b_k^T\u + c_k), \\ 
\nabla_{u}L(\u,\y)/2 & = & \u^T\A_0 +\b_0^T + \sum_{k\in M} y_k (\u^T\A_k + \b_k^T), \\ 
\nabla_{uu}L(\u,\y)/2 & = & \A_0 + \sum_{k\in M} y_k \A_k,
\end{eqnarray*}
for every $(\u,\y) \in \Real^{n + m}$. Hence $L(\u,\y) = \inf\{L(\u',\y) : \u' \in \Real^n\} \equiv \varphi(\y)$ 
if and only if $\nabla_{u}L(\u,\y)/2 = \0^T$ and $\nabla_{uu}L(\u,\y)/2 \in \SymMat^{1+n}_+$ ({\it i.e.},  
$L(\cdot,\y) : \Real^n \rightarrow \Real$ is convex). Since Condition (A) holds if and only if 
$L(\u,\y) = \varphi(\y)$ and and~\eqref{eq:complemetarity} hold as shown in Section 2.1, (A) is equivalent to 
the following condition. \vspace{2mm} \\
{\bf (A') } $(\u,\y) \in \Real^{n+m}$ satisfies 
\begin{eqnarray}
\left.
\begin{array}{l}
\u \in F, \ \y \in G, \ y_kq_k(\u) = 0 \ (k \in M),\\[3pt]  
\displaystyle 
\u^T\A_0+\b_0^T+\sum_{k\in M} y_k (\u^T\A_k + \b_k^T) = \0^T, 
\end{array}
\right\} 
\label{eq:QCQPKKT} 
\end{eqnarray}
(the {\it Karush-Kuhn-Tucker (KKT) condition}), 
and  
\begin{eqnarray}
\A_0+\sum_{k\in M} y_k \A_k & \in & \SymMat^n_+. \label{eq:QCQPglobal}
\end{eqnarray}

\medskip

\noindent
The equivalence of Condition (A) and (A') are well-known \cite{BAZARAA2006}.
(A') is called positive semidefinite condition in \cite{LU2011}. 

\subsection{SDP relaxation of QCQP~\eqref{eq:QCQP00}}

We need to reformulate QCQP~\eqref{eq:QCQP00} to introduce its SDP relaxation. Let 
\begin{eqnarray}
\left.
\begin{array}{l}
\Q_k = \begin{pmatrix} c_k & \b_k^T \\ \b_k & \A_k \end{pmatrix}\in\SymMat^{1+n} \ 
(0 \leq k \leq m),\\[3pt]  
\x   =  \begin{pmatrix} x_0 \\ \u \end{pmatrix} \in \Real^{1+n}, \ 
\H \ = \ \begin{pmatrix} 1 & \0^T \\ \0 & \O \end{pmatrix} \in  \SymMat^{1+n}.
\end{array} 
\right\} \label{eq:notation2}
\end{eqnarray}
Then 
\begin{eqnarray}
q_k(\u) = \Q_k \bullet \begin{pmatrix} 1 \\ \u \end{pmatrix}\begin{pmatrix} 1 \\ \u \end{pmatrix}^T 
\mbox{for every } \u \in \Real^n \ (0 \leq k \leq m), \label{eq:equivalence1}
\end{eqnarray}
and we can rewrite QCQP~\eqref{eq:QCQP00} as 
\begin{eqnarray}
\tilde{\zeta} & = & \inf \left\{ \Q_0 \bullet \x\x^T : \x\x^T \in \widetilde{F} \right\}. \label{eq:QCQP}
\end{eqnarray}
Here 
\begin{eqnarray*}
\widetilde{F} & = & \left\{ \X \in \SymMat^{1+n} :  
\Q_k \bullet \X \leq 0 \ (k\in M),%  \\[3pt] 
\ \H \bullet \X = 1  
\right\}.  
\end{eqnarray*}
We notice that the equality constraint $\H \bullet \x\x^T = 1$ does not 
specify $x_0 = +1$, instead, it requires either $x_0 = +1$ or $x_0 = -1$. 
We see that  if $\x$ is a feasible solution of QCQP~\eqref{eq:QCQP} with the objective value 
$\Q_0 \bullet \x\x^T$, then $-\x$ is a feasible solution of QCQP~\eqref{eq:QCQP} with 
the same objective value $\Q_0 \bullet (-\x)(-\x)^T = \Q_0 \bullet \x\x^T$. 
Thus,   
 the constraint $x_0 \geq 0$ can be implicitly added to QCQP~\eqref{eq:QCQP}.

By replacing $\x\x^T$ by a matrix variable $\X \in \SymMat^{1+n}_+$, we obtain an SDP 
relaxation of QCQP~\eqref{eq:QCQP} and its dual:
\begin{eqnarray}
\eta_p & = & \inf \left\{\Q_0 \bullet \X : 
\begin{array}{ll} 
\X \in \SymMat^{1+n}_{+}, \ \X \in \widetilde{F} 
\end{array}
\right\}.  \label{eq:PSDP} \\ [3pt]
\eta_d & = & \sup \left\{ s : 
\begin{array}{ll} 
\displaystyle \S(\y,s) \equiv \Q_0+\sum_{k \in M}y_k\Q_k - s \H \in \SymMat^{1+n}_+, 
\ \y \in G 
\end{array}
\right\}. \label{eq:DSDP}
\end{eqnarray}
If we add the constraint that rank$\X = 1$ or equivalently $\X = \x\x^T$, 
then the primal SDP~\eqref{eq:PSDP} is equivalent to QCQP~\eqref{eq:QCQP} (hence equivalent 
to~\eqref{eq:QCQP00}). 

For every feasible solution $\X$ of \eqref{eq:PSDP} and every feasible solution $(\y,s)$ of \eqref{eq:DSDP}, 
we observe that 
\begin{eqnarray*}
0 & \leq & \S(\y,s) \bullet \X \\
& = & (\Q_0+\sum_{k \in M}y_k\Q_k - s \H) \bullet \X \\ 
& = & \Q_0\bullet\X +\sum_{k \in M}y_k\Q_k\bullet \X - s \ \leq \ \Q_0\bullet\X - s.  
\end{eqnarray*}
Hence $\eta_d \leq \eta_p$, and the condition 
\begin{eqnarray}
\left. 
\begin{array}{l}
\X \in \SymMat^{1+n}_+, \ \X \in \widetilde{F}, \ \y \in G, \ \ y_k(\Q_k\bullet \X) = 0 
\ (k \in M), \\[3pt]
\S(\y,s) \in \SymMat^{1+n}_+, \  \S(\y,s) \bullet \X = \O. 
\end{array}
\right\} \label{eq:SDPKKT}
\end{eqnarray}
({\em the Karush-Kuhn-Tucker (KKT) condition} for the primal SDP~\eqref{eq:PSDP}) is equivalent to 
\begin{eqnarray}
\left.
\begin{array}{l}
\mbox{$\X$ is an optimal solution of SDP~\eqref{eq:PSDP},} \\[2pt]
\mbox{$(\y,s)$ is optimal solutions of SDP~\eqref{eq:DSDP}, and $\eta_d = \eta_p$}.
\end{array}
\right\} \label{eq:PDoptsolutions}
\end{eqnarray}
Therefore, we can rewrite Condition (B) as \vspace{2mm} \\
{\bf (B') } $\exists (\x,\y,s)$; $\X = \x\x^T$ and $(\y,s)$ satisfy \eqref{eq:SDPKKT}. 

\bigskip

The following result is well-known \cite{FUJIE1997,SHOR1987,SHOR1990}. 
% \vspace{2mm}\\ 
\prop \label{prop:SDP} { \ } \vspace{-2mm} 
Let 
$ % \begin{eqnarray*}
\widetilde{M}^- = \left\{ k \in M : \X \in \widetilde{F} \ \mbox{and } 
\Q_k \bullet \X < 0 \ \mbox{for some } \X \in \SymMat^{1+n}_+ \right\}.
$ % \end{eqnarray*}
If 
\begin{eqnarray}
\left\{  \X \in \widetilde{F}: \X \in \SymMat^{1+n}_{++}, \ \Q_k \bullet \X < 0 \ (k \in \widetilde{M}^-) \right\} 
\not= \emptyset. \label{eq:SlatrCondition} 
\end{eqnarray}
({\rm the (generalized) Slater constraint qualification}) 
holds and $\X$ is an optimal solution of ~\eqref{eq:DSDP}, then 
there exists a $(\y,s) \in \Real^{m+1}$ such that $(\X,\y,s)$ 
satisfies~\eqref{eq:SDPKKT}.   
\eprop

%
%\input sect3.tex
%!TEX root = ./main.tex
%\section{Equivalence of Conditions (A), (B) and (C)}
\section{Proofs of the relations in \eqref{eq:relations}}.

Proof of (A) $\Leftrightarrow$ (B) is given in Section 3.1. 
(B) $\Rightarrow$ (C), (C) $\Rightarrow$ (D) and the relation that 
(D) $\Leftarrow$ (E) if (B) holds are obvious. 
(D) $\Leftarrow$ (E) also follows directly from $\eta_d \leq \eta_p \leq \zeta$. 
By Proposition~\ref{prop:SDP}, we see that (B) holds if  
the Slater constraint qualification~\eqref{eq:SlatrCondition} and (C) hold.   
%It is well-known that 
The relation `(C) $\Leftarrow$ (D) if QCQP has a solution' and the equivalence of   (E) and (F) are well-known, but their proofs are presented %we prove
 in Section 3.2 %The equivalence of (E) and (F) is well-known, but we 
and  Section 3.3, respectively, for completeness.  

\subsection{Proof of (A) $\Leftrightarrow$ (B) and a related result}

We have already seen the equivalence of (A) and (A') and the equivalence of 
(B) and (B') in Section 2. 
Hence, it suffices to show the equivalence of (A') and (B'). 
Take an arbitrary $(\u,\y) \in \Real^{n+m}$. 
By the relation~\eqref{eq:equivalence1}, we see that 
\begin{eqnarray}
\left. 
\begin{array}{l}
\u \in F, \ \y \in G, \\[3pt] 
y_kq_k(\u) = 0 \ (k \in M)
\end{array}
\right\}
\Leftrightarrow
\left\{ 
\begin{array}{l}
\X = \begin{pmatrix} 1 \\ \u \end{pmatrix}\begin{pmatrix} 1 \\ \u \end{pmatrix}^T  \in \widetilde{F}, 
\ \y \in G, \\[3pt]
y_k(\Q_k\bullet \X) = 0 \ (k \in M). 
\end{array}
\right. \label{eq:relation1} 
\end{eqnarray}
It remains to show that 
\begin{eqnarray}
\left.
\begin{array}{l}
\nabla_{u}L(\u,\y) = \0^T, \ \displaystyle \A_0+\sum_{k\in M}y_k\A_k \in \SymMat^{n}_+,\\ 
\displaystyle s  = \big(\b_0^T+\sum_{k\in M}y_k\b_k^T\big)\u+\sum_{k \in M}y_kc_k.
\end{array}
\right\}
\Leftrightarrow  
\left\{
\begin{array}{l}
\S(\y,s) \in \SymMat^{1+n}_+, \\
\S(\y,s) \bullet \begin{pmatrix} 1 \\ \u \end{pmatrix}\begin{pmatrix} 1 \\ \u \end{pmatrix}^T = 0,  
 \end{array}
 \right.\label{eq:relation2}
\end{eqnarray}
which can be proved from the following relations: 
\begin{eqnarray*}
\S(\y,s) & = & \Q_0+\displaystyle \sum_{k\in M} y_k \Q_k - s \H  \\ 
%& \LBlue{=} & \LBlue{
%\begin{pmatrix} c_0 & \b_0^T \\ \b_0 & \A_0 \end{pmatrix} 
%\displaystyle + \sum_{k\in M} y_k\begin{pmatrix} c_k & \b_k^T \\ \b_k & \A_k \end{pmatrix} 
%\displaystyle - s \begin{pmatrix} 1 & \0^T \\ \0 & \O \end{pmatrix}
%} \\ 
& = & \begin{pmatrix} \displaystyle \sum_{k\in M}y_kc_k -s 
& \displaystyle \b_0^T +\sum_{k\in M}y_k\b_k^T \\ 
\displaystyle \b_0 +\sum_{k\in M}y_k\b_k & \displaystyle \A_0 +\sum_{k\in M}y_k\A_k 
\end{pmatrix}, 
\end{eqnarray*}
\begin{eqnarray*}
\S(\y,s) \bullet \begin{pmatrix} 1 \\ \u \end{pmatrix}\begin{pmatrix} 1 \\ \u \end{pmatrix}^T & = & 0 
\Leftrightarrow \S(\y,s) \begin{pmatrix} 1 \\ \u \end{pmatrix} = \0 \ \mbox{if } 
\S(\y,s) \in \SymMat^{1+n}_+, 
\end{eqnarray*}
\begin{eqnarray*}
\S(\y,s)\begin{pmatrix} 1 \\ \u \end{pmatrix} 
& = &
\begin{pmatrix} \displaystyle \sum_{k\in M}y_kc_k -s 
& \displaystyle \b_0^T+\sum_{k\in M}y_k\b_k^T \\ 
\displaystyle \b_0+\sum_{k\in M}y_k\b_k & \displaystyle \A_0+\sum_{k\in M}y_k\A_k 
\end{pmatrix}
\begin{pmatrix} 1 \\ \u \end{pmatrix}\\
%& \LBlue{ = } & 
%\LBlue{
%\begin{pmatrix}\displaystyle -s+c_0 +\sum_{k\in M}y_kc_k + \b_0^T\u+\sum_{k\in M}y_k\b_k^T\u \\ 
%\displaystyle \b_0+\sum_{k\in M}y_k\b_k+\A_0\u+\sum_{k\in M}y_k\A_k\u \end{pmatrix}}\\
& =&\begin{pmatrix}
\displaystyle-s+
% \sum_{k\in M}y_kc_k+\b_0^T\u+\sum_{k\in M}y_k\b_k^T\u \\ 
\big(\b_0^T+\sum_{k\in M}y_k\b_k^T\big)\u+\sum_{k \in M}y_kc_k \\
\nabla_{u}L(\u,\y)^T/2 \end{pmatrix}, 
\end{eqnarray*}
\begin{eqnarray*}
\begin{pmatrix} 1 & \u^T \\ \0 & \I \end{pmatrix} \S(\y,s)\begin{pmatrix} 1 & \0^T \\ \u & \I\end{pmatrix}
& = & \begin{pmatrix} 1 & \u^T \\ \0 & \I \end{pmatrix}
\begin{pmatrix}\displaystyle \sum_{k\in M}y_kc_k-s&\displaystyle\b_0^T+\sum_{k\in M}y_k\b_k^T \\ 
\displaystyle \b_0 +\sum_{k\in M}y_k\b_k & \displaystyle \A_0+\sum_{k\in M}y_k\A_k   
\end{pmatrix} 
\begin{pmatrix} 1 & \0^T \\ \u & \I \end{pmatrix} \nonumber \\ [3pt]
%& \LBlue{=} & 
%\LBlue{
%\begin{pmatrix} 1 & \u^T \\ \0 & \I \end{pmatrix}
%\begin{pmatrix} \displaystyle -s + c_0+\sum_{k\in M}y_kc_k  + \b_0^T\u+\sum_{k\in M}y_k\b_k^T\u 
%& \displaystyle \b_0^T+\sum_{k\in M}y_k\b_k^T \\ 
%\nabla_{u}L(\u,\y)^T/2
%& \displaystyle \A_0+\sum_{k\in M}y_k\A_k 
%\end{pmatrix}} \nonumber \\ [3pt]
& = & 
\begin{pmatrix} 
\displaystyle -s + \big(\b_0^T+\sum_{k\in M}y_k\b_k^T\big)\u+\sum_{k \in M}y_kc_k
% + \sum_{k\in M}y_kc_k + \b_0^T\u+\sum_{k\in M}y_k\b_k^T\u 
% & \displaystyle \b_0^T+\sum_{k=1}^my_k\b_k^T + \u^T\A_0+\sum_{k=1}^my_k\u^T\A_k \\ 
& \nabla_{u}L(\u,\y)/2 \\
% \displaystyle \b_0+\sum_{k=1}^my_k\b_k + \A_0\u+\sum_{k=1}^my_k\A_k\u 
\nabla_{u}L(\u,\y)^T/2 
& \displaystyle \A_0+\sum_{k\in M}y_k\A_k 
\end{pmatrix}. %\label{eq:SMat}
\end{eqnarray*}
\qed

We now consider the following two sufficient conditions for Conditions (A') and (B'), 
respectively. \vspace{2mm}\\
{\bf ($\overline{\rm A}$) } 
$\exists (\u,\y) \in \Real^{n+m}$;  \eqref{eq:QCQPKKT} 
and  
$ \displaystyle \A_0+\sum_{k\in M} y_k \A_k  \in  \SymMat^n_{++}$ hold. \\ % \vspace{2mm} \\
{\bf ($\overline{\rm B}$) } $\exists (\X,\y,s)$; \eqref{eq:SDPKKT} and rank$\S(\y,s) = n$ hold. 
\vspace{2mm} \\
Condition ($\overline{\rm A}$) implies that $\u$ is the unique global minimizer of QCQP~\eqref{eq:QCQP00}, 
while ($\overline{\rm B}$) has been used to identify a class of QCQPs whose SDP relaxation is exact 
in the papers~\cite{AZUMA2022a,AZUMA2022}. 
These two conditions are equivalent. In fact, 
the proof of (A') $\Leftrightarrow$ (B') above can be  modified in a straightforward manner to show %for 
the equivalence relation 
\begin{eqnarray*}
\left.
\begin{array}{l}
\nabla_{u}L(\u,\y) = \0^T, \ 
\displaystyle \A_0 +\sum_{k\in M}y_k\A_k \in \SymMat^{n}_{++},\\ 
\displaystyle s  = \big(\b_0^T+\sum_{k\in M}y_k\b_k^T\big)\u + \sum_{k\in M}y_kc_k,
\end{array}
\right\}
\Leftrightarrow  
\left\{
\begin{array}{l}
\S(\y,s) \in \SymMat^{1+n}_+, \ \mbox{rank}\S(\y,s) = n, \\
\S(\y,s) \bullet \begin{pmatrix} 1 \\ \u \end{pmatrix}\begin{pmatrix} 1 \\ \u \end{pmatrix}^T = 0,  
 \end{array}
 \right.
\end{eqnarray*}
which together with~\eqref{eq:relation1} implies the desired result.  
\qed

\subsection{Proof of `(C) $\Leftarrow$ (D) if QCQP~\eqref{eq:QCQP00} has a minimizer'}

Assume that $\eta_p = \zeta$ and QCQP~\eqref{eq:QCQP00} has a minimizer $\u \in \Real^n$. 
Then $\X = \begin{pmatrix} 1 \\ \u\end{pmatrix}\begin{pmatrix} 1 \\ \u\end{pmatrix}^T$ is a feasible 
solution of the primal SDP~\eqref{eq:PSDP} with the objective value $\zeta = \eta$. Hence $\X$ is a 
rank-$1$ optimal solution of the primal SDP~\eqref{eq:PSDP}. Therefore, (C) holds. 
\qed

\subsection{Proof of (E) $\Leftrightarrow$ (F)} 

For the inner minimization of the Lagrangian dual~\eqref{eq:LagrangianDual}, we observe that 
\begin{eqnarray*}
\varphi(\y) & = & \inf\left\{ L(\u,\y) : \u \in \Real^n \right\} \\ 
            & = & \inf_{\u \in \Real^n} \left\{ L(\u,\y) : 
            \begin{array}{l}
            \displaystyle
            \0^T = L_{u}(\u,\y)/2 \equiv \u^T\big( \A_0 + \sum_{k \in M} y_k\A_k \big) + 
            \b_0^T + \sum_{k \in M} y_k\b_k^T, \\ 
            \displaystyle
            \SymMat^n_+ \ni L_{uu}(\u,\y)/2 \equiv \A_0 + \sum_{k \in M} y_k\A_k
            \end{array}          
             \right\}\\
             & = & \inf_{\u \in \Real^n}\left\{ \big(\b_0 + \sum_{k\in M}y_k\b_k\big)^T\u + \sum_{k\in M}y_kc_k : 
            \begin{array}{l}
            \displaystyle
            \u^T\big( \A_0 + \sum_{k \in M} y_k\A_k \big) + 
            \b_0^T + \sum_{k \in M} y_k\b_k^T = \0^T, \\ 
            \displaystyle
            \A_0 + \sum_{k \in M} y_k\A_k \in \SymMat^n_+ 
            \end{array}          
             \right\}\\
             & = & \inf_{\u \in \Real^n}\left\{ s : 
            \displaystyle
            \S(\y,s) \in \SymMat^{n+1}_+, \ 
            \S(\y,s) \bullet \begin{pmatrix}1 \\ \u \end{pmatrix}\begin{pmatrix}1 \\ \u \end{pmatrix}^T        
             \right\} \ \mbox{(by~\eqref{eq:relation2})}. 
\end{eqnarray*}
Therefore, % we obtain 
\begin{eqnarray*}
	\varphi & = & \sup_{\y \in G} \varphi(\y) = \sup \left\{ s : \S(\y,s) \in \SymMat^{n+1}_+, \ \y \in G\right\} = \eta_d.  \qed
\end{eqnarray*} 
%
%\input sect4Examples.tex
%!TEX root = ./main.tex
\section{Examples}

In this section, we present four QCQP examples 
%, Examples~\ref{examp:QCQP41},\ref{examp:QCQP42}, \ref{examp:QCQP43} and \ref{examp:QCQP44},  
to supplement the relations in~\eqref{eq:relations} and the discussions thus far. 
Table 1 summarizes their characteristics. 
%See Table 1 for the summary of these examples.   

\begin{table}[h!]
\begin{center}
\scriptsize{
\begin{tabular}{|c|c|c|c|c|l|c|c|}
\hline
% Equivalent 
Conditions         &     & Opt. sol.  & SDP KKT  & Conditions (D), (E)  & \\

(A') and (B')      &  (C) & $(\y,s)$ of~\eqref{eq:DSDP} & Cond.~\eqref{eq:SDPKKT} & $\eta_d \leq \eta_p \leq \zeta < \infty $  & Example \\
\hline                
$\circ$               & $\circ$              & $\exists$ & $\circ$  & $-\infty< \eta_d=\eta_p = \zeta$  
&Ex.~\ref{examp:QCQP41}: $\alpha\leq 2,3\leq\alpha<4$,$4<\alpha$\\ 
&                     &                      &                      &
&$\exists$QCQP minimizer \\ 
\hline                
$\times$              & $\circ$              & $\not\exists$ & $\times$ & $-\infty< \eta_d= \eta_p = \zeta$ 
&Ex.~\ref{examp:QCQP41}: $\alpha=4$, $\exists$QCQP minimizer \\
\hline                
$\times$              & $\times$             & $\exists$ & $\circ$ & $-\infty< \eta_d=\eta_p < \zeta$ 
&Ex.~\ref{examp:QCQP41}: $2 < \alpha <3$, $\exists$QCQP minimizer \\
\hline                
$\times$              & $\times$             & $\not\exists$ & $\times$ & $-\infty = \eta_d= \eta_p < \zeta$ 
&Ex.~\ref{examp:QCQP42}, $\not\exists$QCQP minimizer \\
\hline                
$\times$              & $\times$             & $\exists$ & $\circ$ & $-\infty< \eta_d=\eta_p = \zeta$ 
&Ex.~\ref{examp:QCQP43}, $\not\exists$QCQP minimizer \\
\hline                
$\times$              & $\circ$             & $\exists$ & $\times$ & $-\infty< \eta_d < \eta_p = \zeta$ 
&Ex.~\ref{examp:QCQP44}, $\exists$QCQP minimizer\\
\hline                
\end{tabular}	
}
\end{center}
\caption{
The characteristics of Examples~\ref{examp:QCQP41},\ref{examp:QCQP42}, \ref{examp:QCQP43} and \ref{examp:QCQP44}.  
Condition (A') is equivalent to (A) $\exists$ a saddle point $(\u,\y)$ of the Lagrangian function.
(B') is equivalent to (B) the primal 
SDP~\eqref{eq:PSDP} is exact, $\exists$ an optimal solution $(\y,s)$ of the dual SDP~\eqref{eq:DSDP} and 
$-\infty < \eta_d = \eta_p = \zeta < \infty$. (C)
the  primal SDP~\eqref{eq:PSDP} is exact. $\eta_d$, $\eta_p$ and $\zeta$ 
denote optimal values of the dual SDP~\eqref{eq:DSDP}, the  primal SDP~\eqref{eq:PSDP} and QCQP~\eqref{eq:QCQP00}, 
respectively.}
\end{table}
%}

\examp \label{examp:QCQP41}
\begin{eqnarray}
\zeta & = & \inf\left\{ q_0(u) \equiv u^2 : 
\begin{array}{l}
q_1(u) \equiv (u -\alpha)(u-4) \leq 0, \\
q_2(\alpha) \equiv -(u-2)(u-3) \leq 0 
\end{array}
\right\}. 
\label{eq:QCQPexample41}
\end{eqnarray}
Here $n =1$ and $m=2$. 
This example illustrates relations among Conditions (A') (equivalent to (A)), (B') (equivalent to (B)), (C) 
and (D), and shows that the Slater constraint 
qualification~\eqref{eq:SlatrCondition} is necessary for (B') $\Leftarrow$ (C).
The Lagrangian function is written as 
\begin{eqnarray*}
%L(\x,\y) & = & x_2^2 - y_1(x_2 -\alpha x_1)(x_2-4x_1) + y_2(x_2-2x_1)(x_2-3x_1) - y_3(x_1^2 -1)\\
%& = & x_2^2 - y_1(x_2^2 - (\alpha+4)x_1x_2 + 4\alpha x_1^2) + y_2(x_2^2 - 5x_1x_2 +6x_1^2) - y_3(x_1^2 -1). 
L(u,\y) & = & q_0(u) + y_1q_1(u) + y_2q_2(u) = 
u^2 + y_1(u -\alpha)(u-4) - y_2(u-2)(u-3). 
\end{eqnarray*}
The KKT condition is written as 
\begin{eqnarray}
%& & 0 \geq  (u -\alpha)(u-4), \ 0  \geq  -(u-2)(u-3), \ 0 \geq y_1, \ \ 0 \geq y_2, \\ 
%& & 0 = y_1(u -\alpha)(u-4), \ 0 = y_2(u-2)(u-3), \\
\left.
\begin{array}{l}
0 \geq q_k(u), \ 0 \leq y_k, \ y_kq_k(u) = 0 \ (k=1,2), \\ [3pt] 
0 = \nabla_{u}L(u,\y) = 2u + y_1(2u-(\alpha+4)) - y_2(2u-5). 
\end{array}
\right\} 
\label{eq:example41KKT}
\end{eqnarray}
The second order sufficient condition~\eqref{eq:QCQPglobal} for  global optimality is written as 
\begin{eqnarray*}
% \nabla_{uu}L(\u,\y)
\nabla_{uu}L(\u,\y) & = & 2+2y_1 -2y_2 \geq 0. 
\end{eqnarray*}
Letting
\begin{eqnarray*}
& & \Q_0 = \begin{pmatrix} 0 & 0 \\ 0 & 1 \end{pmatrix} \in \SymMat^2, \ 
\Q_1 = \begin{pmatrix} 4\alpha & -(4+\alpha)/2 \\ -(4+\alpha)/2 & 1 \end{pmatrix} \in \SymMat^2,\\ 
& & \Q_2 = \begin{pmatrix} -6 & 5/2 \\ 5/2 & -1 \end{pmatrix} \in \SymMat^2, \ 
\H = \begin{pmatrix} 1 & 0 \\ 0 & 0 \end{pmatrix} \in \SymMat^2, 
\end{eqnarray*}
we obtain the SDP relaxation~\eqref{eq:PSDP} and its dual~\eqref{eq:DSDP}. 

For every $\alpha \in \Real$, QCQP~\eqref{eq:QCQPexample41} has a unique global minimizer $u$ 
with the optimal value $\zeta$ such that  
\begin{eqnarray*}
\begin{array}{lllll}
 u = \alpha, & \zeta = \alpha^2 &   \mbox{if $\alpha < 2$}, \\ % & \mbox{ --- Case (i)}, \\ 
 u = 2, & \zeta = 4 &   \mbox{if $\alpha = 2$}, \\ % & \mbox{  --- Case (ii)}, \\ 
 u = 3, & \zeta = 9 &  \mbox{if $ 2 < \alpha < 3$}, \\ % & \mbox{  --- Case (iii)}, \\
 u = 3, & \zeta = 9 &  \mbox{if $\alpha = 3$}, \\ %& \mbox{  --- Case (iv)}, \\
 u = \alpha, & \zeta = \alpha^2 &  \mbox{if $3 < \alpha < 4$},\\ % & \mbox{  --- Case (v)}, \\ 
 u= 4, & \zeta = 16 &  \mbox{if $\alpha = 4$}, \\ % & \mbox{  --- Case (vi)}, \\  
 u = 4, & \zeta = 16 &  \mbox{if $4 < \alpha$}. % & \mbox{  --- Case (vii)}. 
 \end{array}
\end{eqnarray*} 
%Assuming that the KKT condition~\eqref{eq:example41KKT} holds at 
%the global minimizer $u$ 
For each case, we can easily check and/or solve the 
KKT condition~\eqref{eq:example41KKT} for $\y = (y_1,y_2)$ and $\nabla_{uu}L(\u,\y)$. 
Also, it is easy to compute the solution $\X$ of the primal SDP relaxation~\eqref{eq:PSDP} 
%can be computed easily 
for each case. Table 1 and Figure 1 summarize the results. 
Except for two cases $2 < \alpha <3$ and $\alpha = 4$, % Condition (A) as well as 
Conditions (A') and (B') hold.  
In case $2 < \alpha <3$, we see that the  primal 
SDP~\eqref{eq:PSDP} has no rank-$1$ optimal solution, {\it i.e.}, 
(C) does not hold. In this case, neither (A)' nor (B') holds. 
% Condition (A') does not hold and 
% SDP~\eqref{eq:PSDP} has no optimal solution with rank $1$, 
% that $\eta_p =\zeta$ holds  In this exceptional case,  
We note that the KKT condition~\eqref{eq:example41KKT} holds but 
$\nabla_{uu}L(\u,\y) < 0$. 
Figure 1 (b) shows how $0 <  \zeta - \eta_p$ and det$\X$ changes as $\alpha$ increases 
in the interval $(2,3)$, where det$\X = 0$ if and only if rank$\X = 1$ since 
$\X$ is a $(1+1) \times (1+1)$ matrix with $X_{00} = 1$. 
In case $\alpha = 4$, the KKT condition~\eqref{eq:example41KKT} does not hold. 
As Figure 1 (a) shows, the Lagrangian multiplier $y_1 \geq 0$, which exists when $\alpha \not=4$, 
tends to $\infty$ as $\alpha \rightarrow 4$ from below and  above.  In this case, 
rank-$1$ $\X = \begin{pmatrix} 1 & 4 \\ 4 & 16 \end{pmatrix} \in \SymMat^2_+$
is a unique feasible solution of the primal SDP~\eqref{eq:PSDP}, but the Slater 
constraint qualification~\eqref{eq:SlatrCondition} in Proposition~\ref{prop:SDP} does not hold. Hence, 
the existence of $(\y,s) \in \Real^{m+1}$ satisfying \eqref{eq:SDPKKT} is not 
guaranteed. In fact, such a $(\y,s) \in \Real^{m+1}$ does not exist, and the dual 
SDP~\eqref{eq:DSDP} has no optimal solution. Therefore, this case shows that Condition (B) is merely 
sufficient, but not necessary for (C) 
when the Slater constraint qualification~\eqref{eq:SlatrCondition} is not satisfied. 

\begin{table}[htp]
%\scalebox{0.9}{
\scriptsize{
\begin{center}
\caption{Summary of Example~\ref{examp:QCQP41}
} 
\label{table:example41}
\vspace{2mm}
\begin{tabular}{|c|c|c|c|c|c|c|c|c|}
\hline
        &\multicolumn{2}{|c|}{QCQP~\eqref{eq:QCQP00}}&\multicolumn{2}{|c|}{KKT~\eqref{eq:QCQPKKT}} & 
&\multicolumn{2}{|c|}{SDP~\eqref{eq:PSDP}}\\
% \cline{6-8}
Parameter&  $u$   & $\zeta$                &  $y_1$&$y_2$&$\nabla_{uu}L(\u,\y)$ 
& $\X$ &$\eta_p$\\
\hline
$\alpha < 2$ & $\alpha$    & $\alpha^2$             & $\displaystyle y_1=\frac{2\alpha}{4-\alpha} > 0$ & $y_2=0$ 
& $\displaystyle\frac{2\alpha+8}{4 - \alpha} > 0$
& $\displaystyle\begin{pmatrix}1 & \alpha \\\alpha & \alpha^2 \end{pmatrix}$        & $\alpha^2$ \\
\hline
$\alpha = 2$ & $2$    & $4$             & $ 2 \leq y_1 \leq 5$ & $y_2 = -4+2y_1$ 
& $10-2y_1 \geq 0$
& $\displaystyle\begin{pmatrix}1 & 2 \\2 & 4 \end{pmatrix}$        & $4$ \\
\hline
$2 < \alpha < 3$ & $3$    & $9$             & $y_1=0$ & $y_2=6 $
& -10 % $\displaystyle\frac{2\alpha+8}{\alpha-4}$
& $\displaystyle\begin{pmatrix} 1 & \displaystyle\frac{4\alpha-6}{\alpha-1} \\ \displaystyle\frac{4\alpha-6}{\alpha-1} 
& \displaystyle\frac{14\alpha-24}{\alpha-1} \end{pmatrix}$        & $\displaystyle \frac{14\alpha-24}{\alpha-1}$ \\
\hline
$\alpha =3 $ & $3$    & $9$             & $2.5 \leq y_1 \leq 6$ & $y_2=6-y_1$ 
& $-10+4y_1 \geq 0$
& $\begin{pmatrix}1 & 3 \\ 3 & 9 \end{pmatrix}$        & $9$ \\
%             &        &                 &                      & 
%& $\geq 0$ if $-6 \leq y_1 \leq -2.5$ &                &  \\
\hline
$ 3 < \alpha < 4$ & $\alpha$    & $\alpha^2$             & $\displaystyle y_1=\frac{2\alpha}{4-\alpha} > 0$ & $y_2=0$ 
& $\displaystyle\frac{2\alpha+8}{4-\alpha} > 0$
& $\displaystyle\begin{pmatrix}1 & \alpha \\\alpha & \alpha^2 \end{pmatrix}$        & $\alpha^2$ \\
\hline
$\alpha = 4$ & $4$    & $16$             & \multicolumn{3}{|c|}{Not hold}
& $\displaystyle\begin{pmatrix}1 & 4 \\ 4 & 16 \end{pmatrix}$        & $16$ \\
\hline 
$ 4 < \alpha $ & $4$    & $16$           & $\displaystyle y_1=\frac{8}{\alpha -4} > 0$ & $y_2 = 0$
& $\displaystyle\frac{2\alpha}{\alpha-4} > 0$
& $\displaystyle\begin{pmatrix}1 & 4 \\ 4 & 16 \end{pmatrix}$        & $16$ \\
\hline
\end{tabular}
\end{center}
}
%}
\end{table}

%\afterpage{
\begin{figure}[t!]  \vspace{-0.3cm} 
\begin{center}
\mbox{ \ }

\vspace{-50mm} 

\includegraphics[width=0.90\textwidth]{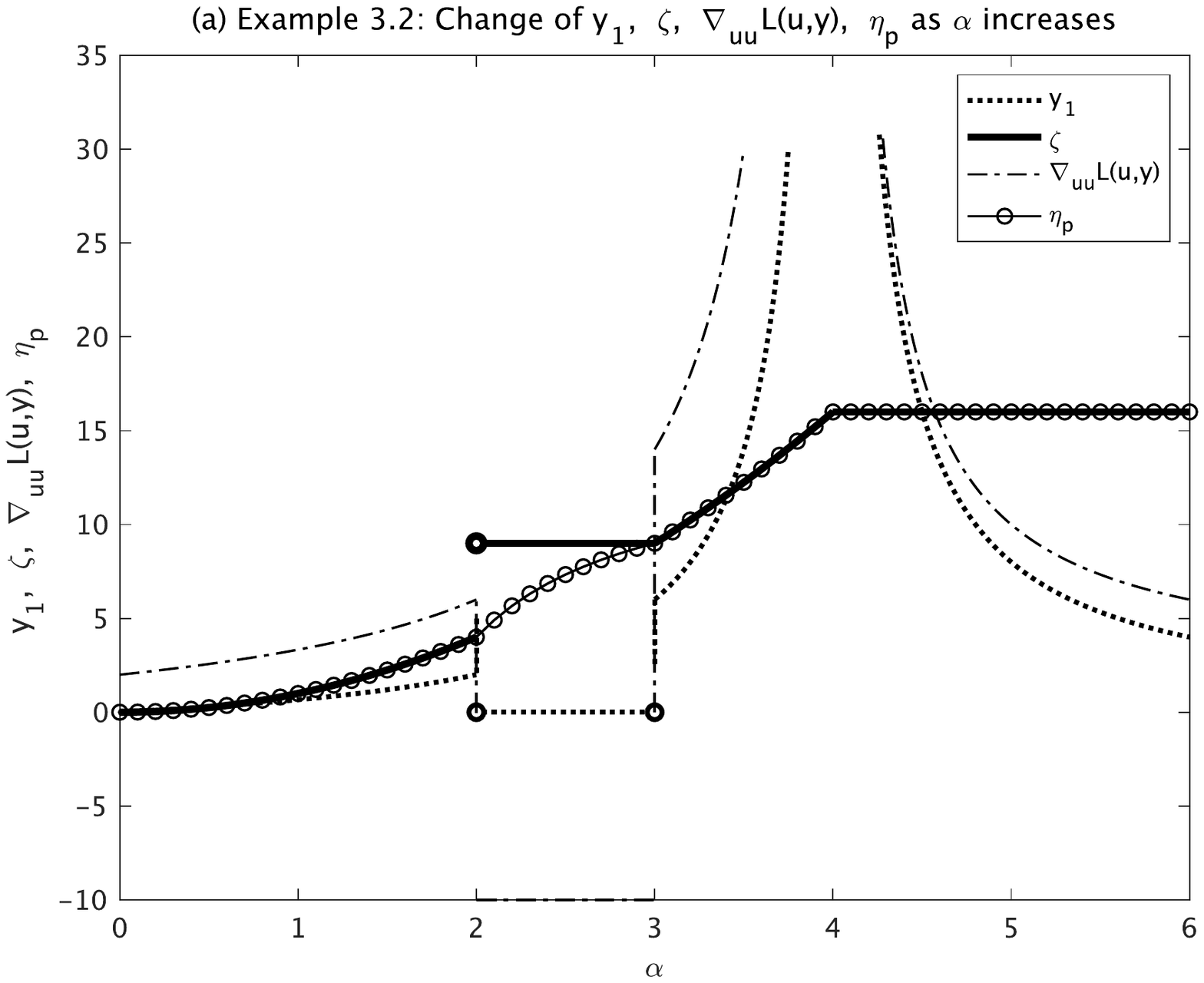}

\vspace{-100mm}

\includegraphics[width=0.90\textwidth]{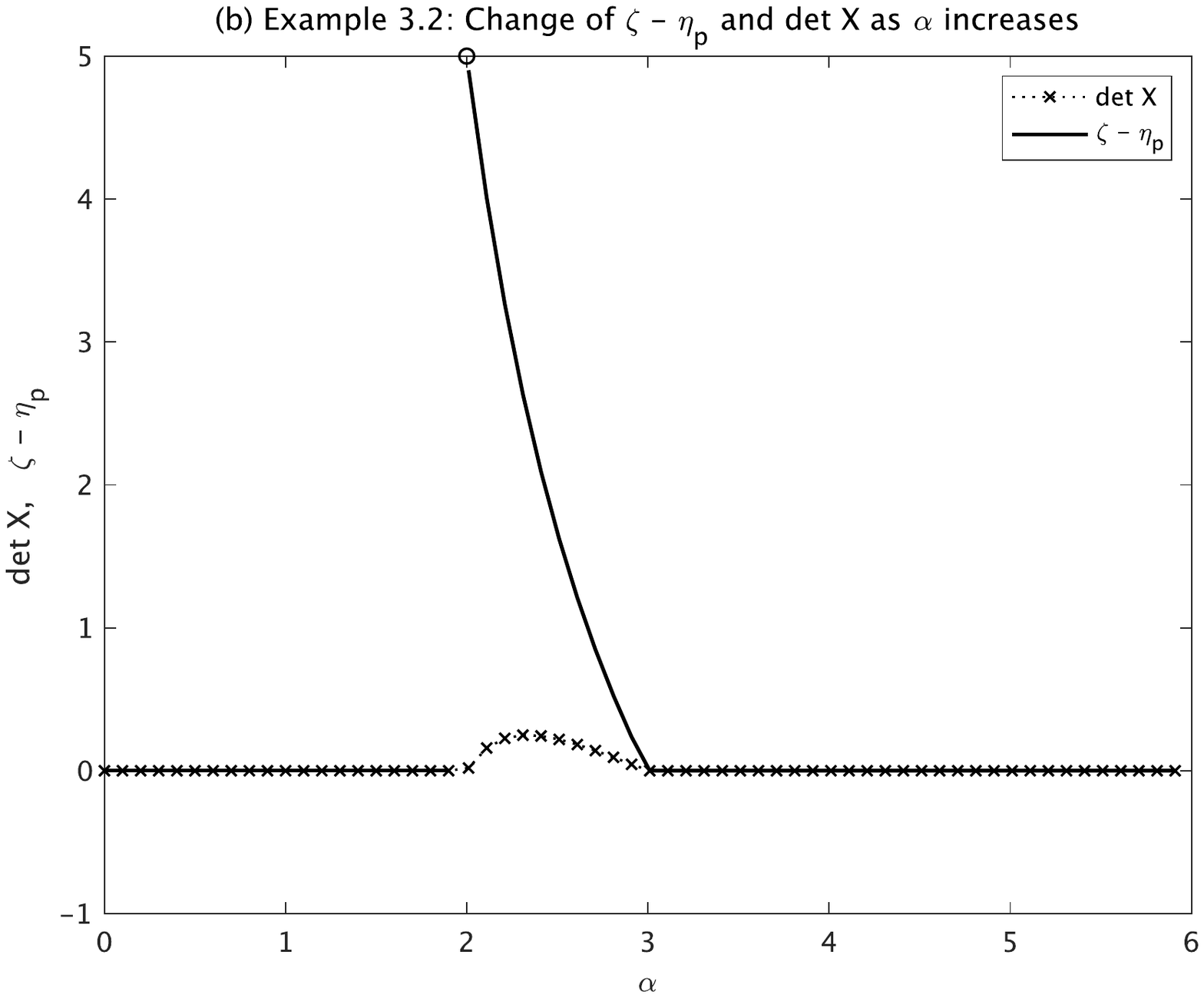}
\end{center}

\vspace{-55mm} 

\caption{
(a) Changes of 
$y_1$, $\nabla_{uu}$, $\zeta$ and $\eta_p$ as $\alpha$ increases from $2$ to $6$ in 
Example 3.2. 
(b) Changes of det$X$ and $\zeta - \eta_p$ as $\alpha$ increases from $2$ to $6$ in 
Example 3.2. Since $X_{11} = 1$, det$X = 0$ if and only if rank$X = 1$. 
}
\end{figure} 
\eexamp

\examp \label{examp:QCQP42}
\begin{eqnarray*}
\zeta & = & 
\inf \left\{ q_0(\u) \equiv 2u_2-2u_1: 
\begin{array}{l}
q_1(\u) \equiv -2u_1 \leq 0, \ q_2(\u) = -2u_2 \leq 0, \\ [2pt]
q_3(\u) \equiv u_1^2-u_2^2 + 1 \leq 0 
\end{array}
\right\}, 
\end{eqnarray*}
where $n=2$ and $m=3$. % are taken. 
This QCQP has a finite optimal value with no minimizer, 
and none of Conditions (A) through (F) hold. 
Obviously, every feasible $(u_1,u_2)$ satisfies  
$u_1 \geq 0$ and $u_2 \geq 1$. 
We also see that $-q_3(\u) = (u_2-u_1)(u_2+u_1) - 1 \geq 0$, which implies that $u_2 > u_1$ and $q_0(\u) > 0$ for every feasible solution $\u$.  %We also see that 
Moreover, 
$\u=(u_1,u_2) = (u_1,\sqrt{u_1^2+1}) \in \Real^2$ 
with $u_1 \geq 0$ is a feasible solution with 
the objective value $q_0(\u) = \sqrt{u_1^2+1}-u_1$, which tends to $0$ as 
$u_1 \rightarrow \infty$. 
Hence $\zeta = 0$ and there is no global minimizer. Letting 
\begin{eqnarray}
\left.
\begin{array}{l}
\Q_0 = \begin{pmatrix} 0 & -1 & 1 \\-1 & 0 & 0 \\ 1 & 0 & 0 \end{pmatrix}, \
\Q_1 = \begin{pmatrix} 0 & -1 & 0 \\-1 & 0 & 0 \\ 0 & 0 & 0 \end{pmatrix},\\[3pt] 
\Q_2 = \begin{pmatrix} 0 & 0 & -1 \\ 0  & 0 & 0 \\ -1 & 0 & 0 \end{pmatrix}, \ 
\Q_3 = \begin{pmatrix} 1 & 0 & 0 \\ 0 & 1 & 0 \\ 0 & 0 & -1 \end{pmatrix}, \
\H = \begin{pmatrix} 1 & 0 & 0 \\ 0 & 0 & 0 \\ 0 & 0 & 0 \end{pmatrix}, 
\end{array}
\right\} 
\label{eq:SDPexample42}
\end{eqnarray}
we obtain the primal SDP relaxation~\eqref{eq:PSDP} and its dual~\eqref{eq:DSDP}. 
We can easily verify that 
if we take $u_1 \geq 0$, then 
\begin{eqnarray*}
\X = \begin{pmatrix} 1 & u_1 & 0 \\ u_1 & u_1^2 & 0 \\ 0 & 0 & u_1^2+1 \end{pmatrix} 
\end{eqnarray*} 
is a feasible solution of the primal SDP~\eqref{eq:PSDP} with the objective 
$\Q_0 \bullet \X = -2u_1$; 
hence $\eta_p = -\infty$. 
Therefore, the dual SDP~\eqref{eq:DSDP} is infeasible. 
In fact, 
the constraints of the dual SDP~\eqref{eq:DSDP} with $\Q_k$ 
$(0 \leq k \leq 3)$ and $\H$ given by~\eqref{eq:SDPexample42} are written as 
\begin{eqnarray*}
& & y_k \geq 0 \ (k=1,2,3), \\ 
& & \S(\y,s) \equiv \Q_0 + \sum_{k=1}^3 y_k\Q_k - s \H % \\
 = \begin{pmatrix} y_3-s        & -\frac{1}{2}-\frac{1}{2}y_1 & \frac{1}{2}-\frac{1}{2}y_2 \\
                 -\frac{1}{2}-\frac{1}{2}y_1  & y_3           & 0             \\
                 \frac{1}{2} -\frac{1}{2}y_2  & 0             & -y_3 
\end{pmatrix} \in \SymMat^3_+. 
\end{eqnarray*}
Clearly, there is no $(\y,s)$ satisfying the constraints. 
\eexamp

\examp \label{examp:QCQP43}
\begin{eqnarray*}
\zeta & = & 
\inf \left\{ q_0(\u) \equiv 2u_2-2u_1: 
\begin{array}{l}
q_1(\u) \equiv -2u_1 \leq 0, \ q_2(\u) = -2u_2 \leq 0, \\ [2pt]
q_3(\u) \equiv u_1^2-u_2^2 + 1 \leq 0, \ q_4(\u) \equiv 2u_1 - 2u_2 \leq 0 
\end{array}
\right\}, 
%\label{eq:QCQPexample42}
\end{eqnarray*}
where $n=2$ and $m=4$. This QCQP is obtained by adding a redundant constraint 
$q_4(\u) \leq 0$ to Example~\ref{examp:QCQP42}, 
so that its optimal value $\zeta$ remains as $0$ and it still has no 
minimizer. However, the characteristics of its SDP relaxation drastically changes. 
Defining 
\begin{eqnarray*}
\Q_4 = \begin{pmatrix} 0 & 1 & -1 \\ 1 & 0 & 0 \\ -1 & 0 & 0 \end{pmatrix} \ 
\end{eqnarray*}
in addition to $\Q_k \ (k=1,2,3)$ and $\H$ given by\eqref{eq:SDPexample42}, 
we obtain the primal SDP relaxation~\eqref{eq:PSDP} and its dual~\eqref{eq:DSDP}. We can easily 
verify that 
\begin{eqnarray*}
\X = \begin{pmatrix} 1 & 0 & 0 \\ 0 & 1 & 0 \\ 0 & 0 & 2 \end{pmatrix} 
\end{eqnarray*} 
is a feasible solution of the primal SDP~\eqref{eq:PSDP} with the objective $\Q_0 \bullet \X = 0$. 
We  see that $(y_1,y_2,y_3,y_4,s) = (0,0,0,1,0)$ is a feasible solution of  the 
dual SDP~\eqref{eq:DSDP} with the objective value $0$. In fact, %we see that
\begin{eqnarray*}
% & & y_k \geq 0 \ (k=1,\ldots,4), \\ 
 \S(\y,s) & \equiv & \Q_0 + \sum_{k=1}^4 y_k\Q_k - s \H  \\
 & = & \begin{pmatrix} y_3-s & -1-y_1+y_4 & 1-y_2-y_4 \\
                 -1-y_1 + y_4 & y_3           & 0             \\
                  1 - y_2 - y_4 & 0           & -y_3 \end{pmatrix} = \O \in \SymMat^3_+.
%\begin{pmatrix}  & -1-y_1+y_4 & 1-y_2-y_4 \\
%                 -1-y_1 + y_4 & y_3           & 0             \\
%                  1 - y_2 - y_4 & 0           & -y_3 
%\end{pmatrix}
%\in \SymMat^3_+. 
\end{eqnarray*} 
Hence, $\eta_d = \eta_p = \zeta = 0$ and the KKT condition~\eqref{eq:SDPKKT} holds. Thus, this example shows 
that even when the the KKT condition~\eqref{eq:SDPKKT} and $\eta_d = \eta_p$ (the strong duality) hold,  
Condition `(D) $\eta_p = \zeta$' does not necessarily 
ensure `Condition (C) primal SDP~\eqref{eq:PSDP} is exact' unless QCQP~\eqref{eq:QCQP00} has a minimizer. 
%In other words,  primal SDP~\eqref{eq:PSDP} may attain the exact optimal value $\zeta$ of QCQP~\eqref{eq:QCQP00} 
%even when it has no rank $1$ optimal solution. 
\eexamp

\examp \label{examp:QCQP44}
\begin{eqnarray*}
\zeta & = & 
\inf \left\{ q_0(\u) \equiv u_3^2: 
\begin{array}{l}
q_1(\u) \equiv u_1^2 \leq 0, \\ [2pt]
q_2(\u) \equiv -2u_1u_2 - 2u_3^2 + 2 \leq 0.
%q_3(\u) \equiv u_1^2-u_2^2 + 1 \leq 0 
\end{array}
\right\},  
%\label{eq:QCQPexample42}
\end{eqnarray*}
where $n=3$ and $m=2$. This example illustrates a case where both Conditions (C) and (D) hold but 
Conditions (E) does not, even when QCQP~\eqref{eq:QCQP00} has a minimizer. 
Obviously, $(u_1,u_2,u_3) = (0,0,1)$ is a minimizer with the objective value $\zeta = 1$. 
Define 
\begin{eqnarray}
\left.
\begin{array}{l}
\Q_0 = \begin{pmatrix} 0&0&0&0\\  0&0&0&0\\ 0&0&0&0\\ 0&0&0&1 \end{pmatrix}, \
\Q_1 = \begin{pmatrix} 0&0&0&0\\  0&1&0&0\\ 0&0&0&0\\ 0&0&0&0 \end{pmatrix},\\ [3pt] 
\Q_2 = \begin{pmatrix} 2&0&0&0\\ 0&0&-1&0\\ 0&-1&0&0\\ 0&0&0&-2 \end{pmatrix}, \ 
\H   = \begin{pmatrix} 1&0&0&0\\  0&0&0&0\\ 0&0&0&0\\ 0&0&0&0 \end{pmatrix}.
\end{array}
\right\} 
\label{eq:SDPexample44}
\end{eqnarray}
Then, we obtain the primal SDP relaxation~\eqref{eq:PSDP} and its dual~\eqref{eq:DSDP}. 
It is easy to see that 
\begin{eqnarray*}
\X = \begin{pmatrix} 1&0&0&1\\  0&0&0&0\\ 0&0&0&0\\ 1&0&0&1 \end{pmatrix}
\end{eqnarray*}
is a rank-$1$ feasible solution of the primal SDP~\eqref{eq:PSDP} with the objective value $\eta_p = 1 = \zeta$; 
hence $\X$ is an optimal solution of~\eqref{eq:PSDP}.   
On the other hand, the constraint of the dual SDP~\eqref{eq:DSDP} 
\begin{eqnarray*}
& &  y_1 \geq 0, \ y_2 \geq 0, \\
& &  \S(\y,s) \equiv \Q_0 + y_1\Q_1 + y_2\Q_2 - s\H \equiv 
\begin{pmatrix} 2y_2-s&0&0&0\\  0&y_1&-y_2&0\\ 0&-y_2&0&0\\ 0&0&0&1-2y_2 \end{pmatrix} \in \SymMat^{1+3}_+
\end{eqnarray*}
holds if and only if $y_1 \geq 0, y_2 = 0$ and $s \leq 0$. Hence, $(y_1,y_2,s) = (0,0,0)$ is an optimal 
solution of the dual SDP~\eqref{eq:DSDP} with the optimal value $\eta_d = 0$. Thus, $0 = \eta_d < 1 = \eta_p = \zeta$ holds. 
\eexamp

%
%\input sect5.tex
%!TEX root = ./main.tex
\section{An extension of Conditions (B), (C), (D) and (E) to doubly nonnegative (DNN) relaxation}

The SDP relaxation  has played a major role in the discussion of Conditions (B), (C), 
(D) and (E).
For QCQP in nonnegative variables, we can strengthen 
those conditions by replacing the SDP relaxation  with a DNN relaxation.
A lower bound provided by the DNN relaxation 
is known to be at least as tight as one by the SDP relaxation in theory 
and is often tighter in practice \cite{KIM2013,KIM2021}.

To discuss conditions for the DNN relaxation, %effectively, we need represent 
 QCQP in nonnegative variables should be  first described, for instance,
%We could rewrite 
by rewriting
QCQP~\eqref{eq:QCQP00}  as 
\begin{eqnarray*}
\zeta & = & \inf\left\{ q_0(\u^+-\u^-) : 
   \begin{array}{l}
   q_k(\u^+-\u^-) + v_k = 0 \ (k \in M), \ v_k \geq 0 \ (k \in M),\\[2pt] 
   \u^+ \geq \0, \ \u^- \geq 0, \ u^+_iu^-_i = 0 \ (1 \leq i \leq n)
   \end{array}
   \right\}, 
\end{eqnarray*}
where $\u^+ \in \Real^n$, $\u^- \in \Real^n$ and $\v \in \Real^m$ are nonnegative variables.
From the transformed QCQP above, a DNN relaxation can be derived.
%and then we could apply DNN relaxation to the transformed QCQP above. But then the description of 
The description of the resulting DNN relaxation would be very complicated.
%too complicated. 
For simplicity of discussion, we instead consider a standard equality form QCQP:
\begin{eqnarray}
\hat{\zeta} & = & \inf\left\{ q_0(\u) : q_k(\u) = 0 \ (k \in M), \ \u \geq \0 \right\}. 
\label{eq:QCQP01} 
\end{eqnarray}
Here $q_k : \Real^n \rightarrow \Real$ $(0 \leq k \leq m)$ denote quadratic functions as used thus far.
% See~\eqref{eq:notation}. 

Introducing redundant quadratic inequalities $u_iu_j \geq 0$ $(1 \leq i,j\leq n)$, 
which can be represented  %we represent 
as a matrix inequality $\begin{pmatrix} x_0^2 & \u^T \\ \u & \u\u^T \end{pmatrix} \geq \O$ with $x_0^2=1$,
and using the notation and symbols given in~\eqref{eq:notation2},   
we first transform QCQP~\eqref{eq:QCQP01} to the following QCQP:
\begin{eqnarray}
\hat{\zeta} & = & 
\inf\left\{ 
\Q_0 \bullet \x\x^T : \x\x^T \geq \O, \ \Q_k \bullet \x\x^T = 0\ (k \in M), \ \H\bullet\x\x^T = 1 
\right\}\nonumber \\ 
& = & \inf\left\{ 
\Q_0 \bullet \x\x^T : \x\x^T \geq \O, \ \x\x^T \in \widehat{F} \right\}. 
\label{eq:QCQP02} 
\end{eqnarray}
Here $\widehat{F} = \left\{ \X \in \SymMat^{1+n} : \Q_k \bullet \x\x^T = 0\ (k \in M), 
\ \H\bullet\x\x^T = 1 \right\}$.
Replacing $\x\x^T$ with a matrix variable $\X \in \SymMat^{1+n}_+$, we now obtain the 
primal SDP relaxation of QCQP\eqref{eq:QCQP02} and its dual as follows. 
\begin{eqnarray}
\hat{\eta}_p & = & \inf\left\{\Q_0\bullet\X:\X\in \SymMat^{1+n}_+\cap\coneN, \ \X\in\widehat{\F} \right\},
\label{eq:PDNN} \\[2pt] 
\hat{\eta}_d & = & \sup\left\{s : \Q_0-\sum_{k \in M}y_k\Q_k-s\H \in \SymMat^{1+n}_+ + \coneN, 
\ (\y,s) \in \Real^{m+1} \right\}, \label{eq:DDNN} 
\end{eqnarray}
which serve as the primal DNN relaxation of QCQP~\eqref{eq:QCQP01} and its dual, 
respectively. Here 
$\coneN = \left\{ \X \in \SymMat^{1+n} : \X \geq \O\right\}$ (the cone of $(1+n) \times (1+n)$ 
nonnegative symmetric matrices). 
The Slater constraint qualification for the primal DNN~\eqref{eq:PDNN} is written as 
\begin{eqnarray}
\left\{\X \in \widehat{F}:\X\in\SymMat^{1+n}_{++}, \ X_{ij}>0 \ (1\leq i,j \leq n)\right\}\not=\emptyset. 
\label{eq:Slater2} 
\end{eqnarray}

Let $\widehat{G} = \Real^{m+1} \times \coneN$. The Lagrangian function 
$L: \Real^{n+m+1} \times \SymMat^{1+n} \rightarrow \Real$ for QCQP\eqref{eq:QCQP02} is defined by 
\begin{eqnarray*}
\widehat{L}(\x,\y,s,\Z) & = & \Q_0\bullet\x\x^T - \sum_{k \in M} y_k \Q_k\bullet\x\x^T - \Z\bullet\x\x^T - s(\H\bullet\x\x^T-1)  \\ 
& & \mbox{for every } (\x,\y,s,\Z) \in \Real^{n+m+1} \times \SymMat^{1+n}. 
\end{eqnarray*}
Hence, the saddle-point problem is: Find a 
$(\x,\y,s,\Z) \in \Real^{1+n} \times \widehat{G}$ such that 
\begin{eqnarray*}
\sup_{ (\y',s',\Z') \in \widehat{G}} L(\x,\y',s',\Z') = L(\x,\y,s,\Z) = \inf_{\x' \in \Real^n} L(\x',\y,s,\Z), 
\end{eqnarray*}
and the Lagrangian dual is: 
\begin{eqnarray*}
\hat{\varphi} & = & \sup_{ (\y,s,\Z) \in \widehat{G}} \ \inf_{\x \in \Real^n} L(\x,\y,s,\Z). 
\end{eqnarray*}

We are now ready to present the following relations. 
\begin{eqnarray}
\left.
\begin{array}{lcl}
& & \hspace{-2mm} \mbox{($\widehat{\rm A}$)  $\exists$ a saddle point $(\x,\y,s,\Z)$ 
of the Lagrangian function $\widehat{L}$}.\\[2pt]
& & \hspace{-1mm} 
\Updownarrow \\[2pt]
& & \hspace{-2mm} \mbox{($\widehat{\rm B}$) the primal DNN~\eqref{eq:PDNN} is exact}, \\[2pt] 
& & \hspace{6mm} \mbox{$\exists$ an optimal solution $(\y,s,\Z)$ of the dual DNN~\eqref{eq:DDNN} and}\\[2pt]
& & \hspace{6mm} \mbox{$\hat{\eta}_d = \hat{\eta}_p$ (the strong duality)}. \\[2pt]
& & % \hspace{2mm} 
\Downarrow \hspace{1mm} \Uparrow \mbox{if ($\hat{\rm a}$) the 
Slater constraint qualification~\eqref{eq:Slater2} holds.}\\[2pt]
& & \hspace{-2mm} \mbox{($\widehat{\rm C}$) the primal DNN~\eqref{eq:PDNN} is exact.}\\[2pt]
& & % \hspace{2mm} 
\Downarrow \hspace{2mm} \Uparrow \mbox{if ($\hat{\rm b}$) 
QCQP~\eqref{eq:QCQP02} has a minimizer $\x$ when $\hat{\zeta}$ is finite.}\\[2pt]
& & \hspace{-2mm} \mbox{($\widehat{\rm D}$) $-\infty < \hat{\eta}_p = \hat{\zeta}  < \infty$.}\\[2pt]
& \mbox{if ($\widehat{\rm B}$) holds} & % \hspace{-1mm} 
\Downarrow \hspace{3mm} \Uparrow \\[2pt]
& & \hspace{-2mm} \mbox{($\widehat{\rm E}$) $-\infty < \hat{\eta}_d = \hat{\zeta} < \infty$} \Longleftrightarrow 
\mbox{($\widehat{\rm F}$) $-\infty < \hat{\varphi} = \hat{\zeta} < \infty$.} 
\end{array}
\right\}
\nonumber % \label{eq:relations2}
\end{eqnarray}
We can prove these relations similarly as in Section~3. The details are omitted. 
%
%\input sect6.tex
%!TEX root = ./main.tex
\section{Concluding remarks}

When QCQP~\eqref{eq:QCQP00} has a finite optimal value, 
the following two cases, (a) and (b), can be considered: 
(a)  the Slater constraint qualification~\eqref{eq:SlatrCondition} holds,
(b)  QCQP~\eqref{eq:QCQP00} has a minimizer, 
If (a) and (b) are satisfied,
then Conditions (A) through (F) for global optimality of QCQP~\eqref{eq:QCQP00} 
are all equivalent.  It was shown in \cite{DUR2017} that (a) is a generic 
property of conic optimization problems. Also, 
if the feasible region of QCQP~\eqref{eq:QCQP00} is bounded, then (b) holds.
Therefore (a) and (b) may be regarded as moderate assumptions to 
avoid special degenerate cases. 
% In particular, (b) is ensured if the feasible region is bounded. 

For (a), however,
we should be more carful as it may be frequently violated in practice. 
Moreover, judging numerically whether (a) is satisfied or not is not an easy task in practice. 
Many computational 
methods including interior-point methods for solving SDPs assume (a) 
for their convergence 
analysis, and often encounter the numerical difficulty when (a) is not satisfied. 

\subsection*{Statements and Declarations}

The authors declare that they have no competing interests.
%

%\newpage
%
%\input appendix.tex
%\input sect5.tex

%\bibliographystyle{plain}
%\bibliography{./enhFOM}

%%%%%

\end{document}